\documentclass[leqno,12pt]{article}
\usepackage{graphicx}
\usepackage{latexsym}
 4

\renewcommand{\epsilon}{\varepsilon}

\newcommand{\eps}{\varepsilon}
\newcommand{\Dkonv}{\stackrel{\mathcal{D}}{\rightarrow}}

\newcommand{\e}{\mathop{\rm E}\nolimits}
\newcommand{\Deq}{\stackrel{\mathcal{D}}{=}}
\newcommand{\n}{\not\in}

\newcommand{\btb}{\begin{table}}
\newcommand{\etb}{\end{table}}

\usepackage{amssymb}

\newcommand{\lt}{\left}
\newcommand{\rt}{\right}

\def\3{\ss}
\def\er{\mathbb{R}}

\def \R{I \!\! R} %Reelle Zahlen
 0

\newcommand{\bea}{\begin{eqnarray*}}
\newcommand{\eea}{\end{eqnarray*}}
\newcommand{\be}{\begin{eqnarray}}
\newcommand{\ee}{\end{eqnarray}}

\def\3{\ss}

\begin{document}

\title{A martingale-transform goodness-of-fit test for the form of the conditional variance}

\author{Holger Dette \\
Ruhr-Universit\"at Bochum \\
Fakult\"at f\"ur Mathematik \\
44780 Bochum, Germany \\
{\small e-mail: holger.dette@ruhr-uni-bochum.de }\\
\and
Benjamin Hetzler\\
Ruhr-Universit\"at Bochum \\
Fakult\"at f\"ur Mathematik \\
44780 Bochum, Germany \\
{\small email: benjamin.hetzler@ruhr-uni-bochum.de}\\
}

\maketitle

\begin{abstract}
In the common nonparametric regression model the problem of testing for a specific parametric form of the variance function is considered.
Recently Dette and Hetzler (2008) proposed a test statistic, which is based on an empirical process of pseudo residuals. The process converges
weakly to a Gaussian process with a complicated covariance kernel depending on the data generating process. In the present paper we consider a
standardized version of this process and propose a martingale transform to obtain asymptotically distribution free tests for the corresponding
Kolmogorov-Smirnov and Cram\'{e}r-von-Mises functionals. The finite sample properties of the proposed tests are investigated by means of a
simulation study.
\end{abstract}



\section{Introduction}
\def\theequation{1.\arabic{equation}}
\setcounter{equation}{0}

We consider the common nonparametric regression model
\begin{equation}\label{1.1}
Y_{i,n} = m(t_{i,n}) + \sigma (t_{i,n}) \varepsilon(t_{i,n}), \quad i = 1, \ldots , n,
\end{equation}
where $\varepsilon_{1,1}, \dots ,\varepsilon_{n,n}$ with $\varepsilon_{i,n}:=\varepsilon(t_{i,n})$ are assumed to form a triangular array of rowwise
independent random variables with mean 0 and variance 1 and $m$
and $\sigma^2$ denote the unknown regression and variance function, respectively. In the regression model (\ref{1.1}) the quantities $0 \leq
t_{1,n} < t_{2,n}<\dots < t_{n,n} \leq 1$ denote the explanatory variables satisfying
\begin{eqnarray}\label{1.2}
\frac{i}{n+1}=\int_0^{t_{i,n}}f\left(t\right)\,dt,\quad i=1,\ldots,n,
\end{eqnarray}
where $f$ denotes a positive density on the interval $[0,1]$ [see Sacks and Ylvisaker (1970)]. Because additional information on the variance
function such as homoscedasticity can improve the efficiency of the statistical inference, several authors have considered the problem of
testing the hypothesis
\begin{equation}\label{1.3}
H_0 : \sigma^2(t) = \sigma^2 (t,\theta);\quad \forall \, \, t \in [0,1],
\end{equation}
in the nonparametric regression model (\ref{1.1}), where $\{ \sigma^2 (\cdot, \theta) \mid \theta \in \Theta \}$ is a given parametric class of
variance functions and $\Theta \subset \mathbb{R}^d$ denotes a finite dimensional parameter space. Most authors consider linear regression
models [see e.g.\ Bickel (1978), Breusch and Pagan (1979), Cook and Weisberg (1983) among others or Pagan and Pak (1993) for a review]. In the
nonparametric regression model (\ref{1.1}) there exist several papers discussing the problem of testing homoscedasticity [see Dette and Munk
(1998), Zhu, Fujikoshi and Naito (2001), Dette (2002) or Liero (2003)]. Recently Dette, van Keilegom and Neumeyer (2007) proposed a test for
the parametric hypothesis (\ref {1.3}), which is based on the difference of two empirical processes of standardized nonparametric residuals
under the null hypothesis and alternative. Weak convergence of the resulting process is shown and -- because the limit distribution is
complicated and depends on certain features of the data generating process -- the consistency of a smoothed bootstrap procedure is established.
Moreover, although the resulting test has nice theoretical and finite sample properties (in particular, it can detect local alternatives
converging to the null hypothesis at a rate $n^{-1/2}$) the approach requires rather strong assumptions regarding the differentiability of the
variance and regression function. Dette and Hetzler (2008) suggested a procedure, which is, on the one hand, able to detect local alternatives
at a rate $n^{-1/2}$ and requires, on the other hand, minimal assumptions regarding the smoothness of the regression and variance function.
These authors proposed to estimate the process \be \label{1.4} S_t(w) = \int^t_0 \biggl(\sigma^2(x) - \sigma^2 (x, \theta^*)\biggr)
\sqrt{w(x)}f(x)\,dx \ee using pseudo residuals [see Gasser, Sroka and Jennen-Steinmetz (1986) or Hall, Kay and Titterington (1990)], where \be
\label{1.5} \theta^* = \arg\min_{\theta \in \Theta} \int^1_0 \biggl ( \sigma^2 (x) - \sigma^2 (x, \theta) \biggr)^2 f(x)\, dx \ee
 is the parameter corresponding to the best
approximation of the function $\sigma^2$ by the parametric class $\{ \sigma^2 (\cdot, \theta) \mid \theta \in \Theta \}$ and $w$ denotes a
weight function [which was actually chosen as $w\equiv1$ by Dette and Hetzler (2008)]. Under very weak smoothness assumptions on the
regression and variance function they proved weak convergence of the estimated process, say $(\hat S_t(w))_{t \in [0,1]}$, to a Gaussian
process. The Kolmogorov-Smirnov and Cram\'{e}r-von-Mises statistic based on $(\hat S_t(w))_{t \in [0,1]}$ were proposed for testing the
hypothesis (\ref{1.3}). Because the covariance kernel of the limiting process depends on the data generating process in a complicated way, a
bootstrap procedure was applied to obtain the critical values.

It is the purpose of the present paper to construct an asymptotically distribution free test for the parametric form of the variance function
which is on the one hand able to detect local alternatives converging to the null hypotheses at a rate $n^{-1/2}$ and on the other hand
requires minimal smoothness assumptions. For this purpose we consider a standardized version of the process discussed by Dette and Hetzler
(2008), where the weight function is estimated from the data. We apply
 the martingale transform proposed by Khmaladze (1981, 1993) in order to obtain a distribution free limiting process. This transformation has been used successfully by
several authors in goodness-of-fit testing problems for hypotheses regarding the regression function [see Stute, Thies and Zhu (1998),
Khmaladze and Koul (2004) or Koul (2006) among others], but to our best knowledge, it has not been studied in the context of testing hypotheses
regarding the variance function. In Section 2 we briefly review the main features of the empirical process proposed by Dette and Hetzler (2008)
and introduce a standardized version of this process which will be the basis for our test statistic. In Section 3 and 4 we consider the
martingale transform and show that the transformed (and standardized) empirical process is asymptotically distribution free. In Section 5 we
discuss several examples and investigate the finite sample properties of a Cram\'{e}r-von-Mises test based on the martingale transformation,
while some of the more technical details are deferred to an appendix.

\section{The basic process based on pseudo residuals}
\def\theequation{2.\arabic{equation}}
\setcounter{equation}{0}

We assume that the regression function $m$, the variance function $\sigma^2$ in (\ref{1.1}), the design density $f$ and the weight function
$w$ in (\ref{1.4}) are Lipschitz continuous of order $\gamma > \frac {1}{2}$ and that the moments of order 8 of the errors $\varepsilon_{i,n}$
exist and are uniformly bounded. In general, the moments of order $j\geq 3$ of the errors may depend on the explanatory variables $t_{i,n}$,
that is
$$
m_j (t_{i,n}) = E \ [\varepsilon^j_{i,n}], \qquad \qquad j = 3,\ldots,8,
$$
and the functions $m_3$ and $m_4$ are also assumed to be Lipschitz continuous of order $\gamma > \frac {1}{2}$. For the sake of a transparent
presentation we consider at the moment linear hypotheses of the form
\begin{eqnarray}\label{2.1}
H_0:\sigma^2\lt(t\rt)=\sum_{j=1}^d\theta_j\sigma_j^2\lt(t\rt),~~~~~~~~\mbox{ for all }t\in\lt[0,1\rt],
\end{eqnarray}
where $\theta_1, \ldots , \theta_d \in \er$ are unknown parameters and $\sigma^2_1, \ldots , \sigma^2_d$ are given linearly independent
functions satisfying
\begin{eqnarray}\label{2.2}
\sigma^2_j \in \, \mbox{Lip}_\gamma [0,1], ~ \quad j = 1, \ldots , d.
\end{eqnarray}
The general case of testing hypotheses of the form (\ref{1.3}) will be briefly discussed at the end of this section. It is is shown in Dette
and Hetzler (2008) that
 the process defined in (\ref{1.4}) can be consistently estimated
 by
 \begin{eqnarray}\label{2.3}
\hat{S}_t (w)=\hat{B}_t^0 (w)-\hat{B}_t^T (w) \hat{A}^{-1}\hat{C} ,
\end{eqnarray}
where the elements of the matrix $\hat A=\lt(\hat a_{ij}\rt)_{1\leq i,j\leq d}$ and the vector $\hat C=\lt(\hat c_1 ,\ldots,\hat c_d\rt)^T$
 are defined by
\begin{eqnarray} \label{2.4}
\hat{a}_{ij} = \frac{1}{n}\sum_{k=1}^n\sigma_i^2\lt(t_{k,n}\rt)\sigma_j^2\lt(t_{k,n}\rt), \quad 1 \le i,j \le d, \\
\hat{c}_i =\frac{1}{n-r}\sum_{k=r+1}^nR_{k,n}^2 \ \sigma_i^2\lt(t_{k,n}\rt) \label{2.5}, \quad 1 \le i \le d,
\end{eqnarray}
respectively, \be\label{2.5a} \hat{B}_t^0 (w) = \frac{1}{n-r}\sum_{j=r+1}^n1_{\{t_{j,n}\leq t\}} \sqrt{w({t_{j,n}})} R_{j,n}^2 \ee
and ${\hat
B}_t (w) = ({\hat B}^1_t (w), \ldots , {\hat B}^d_t (w))^T$ with
\be \hat{B}_t^i (w) =\frac{1}{n}\sum_{j=1}^n1_{\{t_{j,n}\leq t\}} \sqrt{w(t_{j,n})} \ \sigma_i^2\lt(t_{j,n}\rt),
\quad i = 1, \ldots , d \:.\label{2.6} \ee In (\ref{2.5}) and (\ref{2.5a}) the quantities $R_{j,n}$ denote pseudo residuals defined by
\begin{eqnarray}\label{2.7}
R_{j,n}=\sum_{i=0}^r d_iY_{j-i,n},~~j=r+1,\ldots, n,
\end{eqnarray}
where the vector $(d_0, \ldots , d_r)^T\in\R^{r+1}$ satisfies
\begin{eqnarray}\label{2.8}
\sum_{i=0}^rd_i=0,\;\;\;\sum_{i=0}^r d_i^2=1
\end{eqnarray}
and is called difference sequence of order $r$ [see Gasser, Sroka and Jennen-Steinmetz (1986) or Hall, Kay and Titterington (1990) among others]. The
following result was proved in Dette and Hetzler (2008) and provides the asymptotic properties of the process ${\hat S}_t(w)$ for an increasing
sample size.

\bigskip

{\bf Theorem 2.1.} {\it If the conditions stated at the beginning of this section
 are satisfied, then the process
$\{\sqrt{n} ({\hat S}_t(w) - S_t(w))\}_{t \in [0,1]}$ converges weakly in $D[0,1]$ to a centered Gaussian process with covariance kernel $k(t_1,
t_2)$ given by the non-diagonal elements of the matrix $V_2 \Sigma_{t_1, t_2} V_2^T\in\er^{2 \times 2},$ where the matrices $\Sigma_{t_1, t_2}
\in \er^{(d+2)\times (d+2)}$ and $V_2 \in \er^{2 \times (d+2)}$
 are defined by
\begin{eqnarray}
\Sigma_{t_1,t_2} &=&\left(\begin{array}{cccccc}
v_{11} & v_{12} & w_{11} & \cdots & w_{1d}\\
v_{21} & v_{22} & w_{21} & \cdots & w_{2d}\\
w_{11} & w_{21} & z_{11} & \cdots & z_{1d}\\
\vdots & \vdots & \vdots & \ddots & \vdots \\
w_{1d} & w_{2d} & z_{d1} & \cdots & z_{dd}
\end{array}\right),
\label{2.9} \\
&& \nonumber \\
V_2 &=& \lt(I_2|U\rt),\quad\quad U=-\left(\begin{array}{c}
B_{t_1}^T (w) A^{-1} \\
B_{t_2}^T (w) A^{-1}
\end{array} \right) ,
\label{2.10}
\end{eqnarray}
respectively. The vector $B_t^T (w)$ is defined by
\begin{equation} \label{2.10a}
B_t^T (w) =\lt(\int_0^t\sigma_1^2\lt(x \rt)\sqrt{w(x)} f\lt(x\rt)\,dx, \ldots,\int_0^t\sigma_d^2\lt(x\rt) \sqrt{w(x)} f\lt(x\rt)\,dx\rt),
\end{equation}
the elements of the matrix $A=(a_{ij})_{1\leq i,j\leq d}$ are given by
\be\label{2.10b}
a_{ij}=\int_0^1\sigma_i^2(x)\sigma_j^2(x)f(x)\,dx,~\quad 1\leq i,j\leq d,
\ee
the elements of the matrix in (\ref{2.9}) are defined by
\begin{eqnarray*}
v_{ij}&=&\int_0^1 \tau_r(s)
\sigma^4\lt(s\rt)1_{\lt[0,t_i\wedge t_j\rt)}\lt(s\rt) w(s) f\lt(s\rt)\,ds, ~\quad 1\leq i,j\leq 2,\\
w_{ij}&=&\int_0^1 \tau_r(s)
 \sigma^4\lt(s\rt)\sigma_j^2\lt(s\rt)1_{\lt[0,t_i\rt)}\lt(s\rt) \sqrt{w(s)} f\lt(s\rt)\,ds,~\quad 1\leq i\leq 2, 1\leq j\leq d,\\
 z_{ij}&=&\int_0^1 \tau_r(s) \sigma^4\lt(s\rt)\sigma_i^2\lt(s\rt)\sigma_j^2\lt(s\rt) f\lt(s\rt)\,ds,
~\quad 1\leq i,j\leq d
\end{eqnarray*}
with $\tau_r(s) = m_4\lt(s\rt)-1+4\delta_r$, and the quantity $\delta_r$ is given by}
\begin{equation}\label{2.11}
\delta_r=\sum_{m=1}^r\Big(\sum_{j=0}^{r-m}d_jd_{j+m}\Big)^2.
\end{equation}

\bigskip

Note that the
 null hypothesis (\ref{2.1}) (or more generally the hypothesis (\ref{1.3})) is equivalent to $S_t(w) \equiv 0$ \ $\forall \, \,
 t \in [0,1]$, and consequently rejecting (\ref{2.1}) for large values of the Kolmogorov-Smirnov
or Cram\'{e}r-von-Mises statistic
$$
K_n=\sqrt{n}\sup_{t\in\lt[0,1\rt]} | \hat{S}_t(w)|~,~~ G_n = n\int_0^1 | \hat{S}_t(w)|^2 d F_n(t)
$$
yields a consistent test. Here $F_n(t) = \frac{1}{n} \sum^n_{i=1} 1_{\{ t_{i,n} \le t\}}$ is the empirical distribution function of the design
points. Moreover, it is demonstrated by Dette and Hetzler (2008) that this test can detect local alternatives which converge to the null hypothesis
with a rate $n^{- 1/2}$. Because the limiting distribution depends on certain features of the data generating process, these authors proposed a
bootstrap procedure to calculate the critical values.

 If $\lt(A\lt(t,w\rt)\rt)_{t\in\lt[0,1\rt]}$ denotes the limiting process in Theorem 2.1 it follows from the
Continuous Mapping Theorem [see Pollard (1984)] that
\[K_n \Dkonv \sup_{t\in\lt[0,1\rt]}|A\lt(t,w\rt)|~,~~G_n \Dkonv \int_0^1 | {A(t,w)}|^2 d F(t),
\]
where $F$ denotes the distribution function of the design points. Using the Lipschitz continuity of the regression and variance function, it was shown in the proof of Theorem 2.1 that the process $A_n (t,w) =
\sqrt{n} (\hat S_t(w) - S_t(w))$ exhibits the same asymptotic behaviour as the process
\begin{eqnarray} \label{2.12}
\bar{A}_n(t,w)&=&C_n(t,w)-D_n(t,w),
\end{eqnarray}
where
\begin{eqnarray} \label{2.13}
C_n(t,w) &=&
\frac{\sqrt{n}}{n-r}\sum_{i=r+1}^n1_{\{t_{i,n}\leq t\}} \sqrt{w (t_{i,n})} Z_{i,n}, \\
D_n(t,w)&=&
 {B}_t^T (w) {A}^{-1}\Big(\frac{\sqrt{n}}{n-r}\sum_{i=r+1}^{n}Z_{i,n} \ \sigma_j^2(t_{i,n})\Big)_{j=1}^d,
 \end{eqnarray}
the vector $B_t^T (w) =( B_t^1(w),\ldots ,B_t^d(w))$ and the matrix $A =(a_{ij})_{1\leq i,j\leq d}$ are defined in (\ref{2.10a}) and (\ref{2.10b}),
respectively, and the random variables $Z_{i,n}$ are given by $Z_{i,n} = L^2_{i,n} - \e[L^2_{i,n}]$, with
\be
L_{i,n}=\sum_{j=0}^rd_j\sigma(t_{i-j,n})\eps_{i-j,n}. \ee Because $\{ Z_{i,n} \mid i = 1, \dots, n, \ n \in \mathbb{N} \}$ is a triangular
array of $r$-dependent random variables, it follows observing
\be\label{2.14}
\e\lt[Z_{j,n}^2\rt]+2\sum_{m=1}^r\e\lt[Z_{j,n}Z_{j+m,n}\rt]=\lt(m_4\lt(t_{j,n}\rt)-1+4\delta_r\rt) \sigma^4\lt(t_{j,n}\rt)+O\lt(n^{-\gamma}\rt)
\ee
[see Dette and Hetzler (2008)] that the process $\{ C_n (t,w)\}_{t \in [0,1]} $ converges weakly in D[0,1] to the process $W \circ \psi$, where $W$
denotes a Brownian motion and the function $\psi$ is defined by \be \label{2.15} \psi(t)=\int_0^t\beta(x)w(x) f(x)\,dx \ee with
$$
\beta(x)=(m_4(x)-1+4\delta_r)\sigma^4(x).
$$
Note that the transformation $\psi$ depends on the unknown function $\beta$ which is not known, because it contains the variance and the fourth
moments of the innovations $\varepsilon_{i,n}$. In the following we will use the specific weight function $w(x) = 1/\beta (x)$ for which the
function $\psi$ reduces to $\psi (t) = F(t) = \int^t_0 f (x) dx$ and the process $\{ C_n (t,1/\beta)\}_{t \in [0,1]} $
converges weakly to a Brownian motion $W \circ F$. We assume in a first step that the function $\beta$ is known and
 investigate the martingale transformation of the standardized process
\be
A_n^0(t)=C_n^0(t)-D_n^0(t) \: , \label{2.16}
\ee
where $A_n^0(t)=\bar{A}_n(t,1/\beta)$,
\begin{eqnarray} \label{2.17}
C_n^0(t)&=& C_n(t,1/\beta)=\frac{\sqrt{n}}{n-r}\sum_{i=r+1}^n1_{\{t_{i,n}\leq t\}}Z_{i,n}\beta^{-1/2}(t_{i,n}),\\
\label{2.18} D_n^0(t) &=& D_n(t,1/\beta)=B_t^T(1/\beta)A^{-1}\frac{\sqrt{n}}{n-r}\sum_{i=r+1}^{n}Z_{i,n}g(t_{i,n})
\end{eqnarray}
and $g(x)=(\sigma_1^2(x),\ldots,\sigma_d^2(x))^T$.
In a second step we will estimate the function $\beta$ nonparametrically and consider the corresponding processes standardized by this
estimate. More precisely, we will show that the corresponding martingale transform of the process
\be
\label{proc} \Lambda_n(t) =\sqrt{n}\Big ( \hat S_t (1/ \hat \beta) - S_t
(1/\beta)\Big)
\ee
leads to an asymptotically distribution free test, where $\hat \beta $ is an appropriate estimate of the function $\beta$.

\bigskip

{\bf Remark 2.2.} For the problem of testing a general nonlinear hypothesis of the form (\ref{1.3}) we propose to consider the process
$$
{\hat S}_t(w) = {\hat B}_t^0(w) - \frac{1}{n} \sum^n_{i = 1} 1_{\{t_{i,n}\leq t\}}\sigma^2(t_{i,n}, {\hat \theta}) \sqrt{w(t_{i,n})},
$$
where
$$
{\hat \theta} = \arg\min_{\theta \in \Theta} \frac{1}{n-r} \sum^n_{i=r+1} \Bigl( R^2_{i,n} - \sigma^2(t_{i,n}, \theta)\Bigr)^2
$$
is the least squares estimate of the parameter $\theta^*$ defined by (\ref{1.5}).
%$$\inf_{\theta \in \Theta} \int^1_0 \{ \sigma^2(t) - \sigma^2(t, \theta)\}^2 \sqrt{w(t)}f(t) dt$$
In this case it was shown by Dette and Hetzler (2008) that under assumptions of regularity
the process $\{\sqrt{n}
({\hat S}_t(w) - S_t(w))\}_{t \in [0,1]} $ exhibits the same asymptotic behaviour as described in Theorem 2.1 for the linear
 case, where the functions
$\sigma^2_j$ have to be replaced by
$$ \sigma^2_j (t) = \frac{\partial}{\partial \theta_j} \sigma^2 (t, \theta)\Bigl|_{\theta = \theta^*} ,
\quad j = 1, \ldots , d.
$$
Thus all results presented in the following section can be transferred to the nonlinear case using this identification.

\section{The martingale transform of the process $A^0_n$}
\def\theequation{3.\arabic{equation}}
\setcounter{equation}{0}

It follows by similar arguments as given in Dette and Hetzler (2008) that the process $\{ A^0_n (t) \}_{t \in [0,1]}$ defined by (\ref{2.16})
converges weakly in $D [0,1]$, that is
\begin{equation}\label{2.19}
A_n^0\Dkonv W\circ F-B_t^T (1/\beta) A^{-1} V_0=\tilde{R}^0_{\infty},
\end{equation}
where $W$ is a Brownian motion and $V_0$ denotes a centered normal random variable with mean 0 and covariance matrix
$$
L=\int_0^1 g(x) g^T (x) \beta (x) f (x) dx.
$$
 Because the distribution of the process $\tilde R_\infty^0$ is complicated, we consider in the
following section an operator, which transforms the process $\tilde R_\infty^0$ on the martingale part in its corresponding Doob-Meyer
decomposition. Following Khmaladze and Koul (2004) we define a linear operator $T$ such that
\begin{eqnarray}\label{3.1}
TR^0_{\infty} &\Deq& R^0_{\infty} \: , \\
\label{3.2}T(B^T_t (1/ \beta) A^{-1} V_0) &\equiv &0 \: ,
\end{eqnarray}
where the symbol $\Deq $ denotes equality in distribution and the process $R^0_\infty$ is given by $R^0_\infty = W \circ F$. For this purpose we consider the matrix
\be \label{3.2a}
H(t)=\int_t^1 \beta^{-1} (u)g(u)g^T(u)f(u)\,du \:
\ee
and define for a function $\eta$ its transformation $T\eta$ by
\begin{equation}\label{3.3}
(T\eta)(t)=\eta(t)-\int_0^t \beta^{-1/2} (y) g^T(y)H^{-1}(y)\int_y^1 \beta^{-1/2} (z) g(z)\eta(dz)F(dy),
\end{equation}
where only functions are considered such that the integral on the right hand side of (\ref{3.3}) exists. Note that the matrix $H(x)$ is
non-singular for all $x \in [0,1)$ because the functions $\sigma^2_1, \dots, \sigma_d^2$ are linearly independent; see Achieser (1956). If
$\eta$ is a stochastic process on the interval [0,1], the corresponding integral in (\ref{3.3}) is interpreted as an Ito-integral [see
{\O}ksendal (2003)]. A straightforward calculation shows that
\begin{eqnarray*}
&& T (B^T_t (1/\beta) A^{-1} V_0) = 0 \: , \\
&& \mbox{Cov} (TR_\infty (r), T R_\infty (s)) = F (r \wedge s) \: ,
\end{eqnarray*}
which yields for the process defined on the right hand side of (\ref{2.19})
\be \label{3.4}
T\tilde{R}^0_{\infty}\Deq TR^0_{\infty}\Deq R^0_{\infty}\Deq W\circ F
\ee
 (note that $ R^0_\infty$ is a Gaussian process and that the operator $T$ is
linear). The following theorem shows that a similar property holds in an asymptotic sense for the process $\{ A^0_n (t) \}_{t \in [0,1]}$.

\bigskip

{\bf Theorem 3.1.} {\it If the assumptions stated in Section 2 are satisfied, then the transformed process $\{ T A^0_n (t) \}_{t \in [0,1]}$
converges weakly in $D[0,1]$ to a Brownian motion in time $F$, that is}
\[\{TA_n^0(t)\}_{t\in[0,1]}\Dkonv \{W\circ F(t)\}_{t\in[0,1]}.\]

\medskip

{\bf Proof.}
%For the sake of a transparent notation we omit the index $n$ in this proof and in the following sections,
%whenever the dependence on $n$ will be clear from the context.
%In particular we write $t_j$ and $Z_j$ instead of $t_{j,n}$ and $Z_{j,n},$ respectively.
The assertion of the theorem follows from the statements
\begin{eqnarray} \label{3.6}
 TA_n^0&=&TC_n^0,
 \\
 \{TC_n^0(t)\}_{t\in[0,1]}&\Dkonv &\{W\circ F(t)\}_{t\in[0,1]}. \label{3.7}
\end{eqnarray}
For a proof of (\ref{3.6}) we recall the notation $D^0_n = A^0_n - C^0_n$ in (\ref{2.18}) and obtain by a straightforward calculation from the
definitions (\ref{3.3}) and (\ref{2.18})
\begin{eqnarray} \label{z1} \nonumber
T D^0_n(t) &=& TC_n^0(t) - TA_n^0(t)\\
&=& D_n^0(t)-\int_0^t \beta^{-1/2} (y) g^T(y)H^{-1}(y)\int_y^1 \beta^{-1} (z) g(z)g^T(z)F(dz)F(dy)A^{-1}\nonumber\\
\nonumber &&\times \Big(\frac{\sqrt{n}}{n-r}\sum_{k=r+1}^nZ_{k,n}\ g(t_{k,n})\Big)\\ \nonumber &=& D_n^0(t)-\int_0^t \beta^{-1/2} (y)
g^T(y)H^{-1}(y)H(y) F (dy) A^{-1}\Big(\frac{\sqrt{n}}{n-r}\sum_{k=r+1}^nZ_{k,n} \ g(t_{k,n})\Big)\\
&=& D_n^0(t)-B_t^T (1/\beta) A^{-1}\Big(\frac{\sqrt{n}}{n-r}\sum_{k=r+1}^nZ_{k,n} \ g(t_{k,n})\Big) ~=~ 0.
\end{eqnarray}
The process $T C^0_n$ is a sum of $r$-dependent random variables. Therefore, weak convergence of the finite dimensional distributions and
tightness can be shown using similar arguments as in Dette and Hetzler (2008). Thus the assertion follows showing that the covariance kernel of
the limiting process is given by $F (s \wedge t)$. For the calculation of the asymptotic covariances we use the representation
\begin{eqnarray} \label{3.8}
TC_n^0(t)&=& \frac{\sqrt{n}}{n-r} \sum_{i=r+1}^n C_{i,n}(t),
\end{eqnarray}
where
\begin{eqnarray} \label{3.9}\nonumber
C_{i,n}(t)&=&1_{\{t_{i,n}\leq t\}}Z_{i,n}\beta^{-1/2}(t_{i,n})-\int_0^t \beta^{-1/2} (y) g^T(y)H^{-1}(y)1_{\{t_{i,n}\geq y\}}g(t_{i,n})Z_{i,n}\beta^{-1}(t_{i,n})F(dy)\\
&=& C_{i,n}^{\left(1\right)}(t)-C_{i,n}^{\left(2\right)}(t)
\end{eqnarray}
and the last line defines the random variables $C^{(1)}_{i,n} (t)$ and $C^{(2)}_{i,n} (t)$ in an obvious manner. Observing that
\begin{eqnarray} \label{3.10}
\e[Z_{i,n}^2]+2\sum_{m=1}^r\e[Z_{i,n}Z_{i+m,n}]=\beta(t_{i,n})+O(n^{-\gamma})
\end{eqnarray}
[see (\ref{2.14}) or Dette and Hetzler (2008)], it follows for $r \leq s$
\begin{eqnarray*}
\lefteqn{\e[C_{i,n}^{\left(1\right)}(r)C_{i,n}^{\left(1\right)}(s)]+2\sum_{m=1}^r\e[C_{i,n}^{\left(1\right)}(r)C_{i+m,n}^{\left(1\right)}(s)]}\\
&=& 1_{\{t_{i,n}\leq r\}}E[Z_{i,n}^2]\beta^{-1}(t_{i,n})+ 2\sum_{m=1}^r1_{\{t_{i,n}\leq r\}}\e[Z_{i,n}Z_{i+m,n}]\beta^{-1}(t_{i,n})+o(1) =
1_{\{t_{i,n}\leq r\}}+o(1) \: .
\end{eqnarray*}
This implies
$$
\frac{n}{(n-r)^2}\sum_{i=r+1}^{n-r}\e[C_{i,n}^{\left(1\right)}(r)C_{i,n}^{\left(1\right)}(s)]
+2\sum_{m=1}^r\e[C_{i,n}^{\left(1\right)}(r)C_{i+m,n}^{\left(1\right)}(s)] =F(r)+o(1) \: ,
$$
and similar arguments show
\begin{eqnarray*}
\lefteqn{\frac{n}{(n-r)^2}\sum_{i=r+1}^{n-r}\Big(\e[C_{i,n}^{\left(1\right)}(r)C_{i,n}^{\left(2\right)}(s)]+2\sum_{m=1}^r\e[C_{i,n}^{\left(1\right)}(r)C_{i+m,n}^{\left(2\right)}(s)]\Big)}\nonumber\\
&=& \int_0^r \beta^{-1/2} (y) g^T(y)H^{-1}(y)\int_y^r \beta^{-1/2} (x) g(x)F(dx)F(dy)+o(1),\\
&&\\
\lefteqn{\frac{n}{(n-r)^2}\sum_{i=r+1}^{n-r}\Big(\e[C_{i,n}^{\left(1\right)}(s)C_{i,n}^{\left(2\right)}(r)]+2\sum_{m=1}^r\e[C_{i,n}^{\left(1\right)}(s)C_{i+m,n}^{\left(2\right)}(r)]\Big)}\nonumber\\
&=& \int_0^r \beta^{-1/2} (y) g^T(y)H^{-1}(y)\int_y^s \beta^{-1/2} (x) g(x)F(dx)F(dy)+o(1),\\
&&\\
\lefteqn{\frac{n}{(n-r)^2}\sum_{i=r+1}^{n-r}\Big(\e[C_{i,n}^{\left(2\right)}(r)C_{i,n}^{\left(2\right)}(s)]+2\sum_{m=1}^r\e[C_{i,n}^{\left(2\right)}(r)C_{i+m,n}^{\left(2\right)}(s)]\Big)}\nonumber\\
&=& \int_0^r\int_0^s \beta^{-1/2} (y_1) g^T(y_1)H^{-1}(y_1)H(y_1\vee y_2)H^{-1}(y_2) \beta^{-1/2} (y_2)
 g(y_2)F(dy_2)F(dy_1)+o(1) \: .
\end{eqnarray*}
A combination of these results and an application of Fubini's theorem yield
\begin{eqnarray*}
\e[TC_n^0(r)TC_n^0(s)]
&=& \frac{n}{(n-r)^2}\sum_{i=r+1}^{n-r}\Big(\e[C_{i,n}(r)C_{i,n}(s)]+2\sum_{m=1}^r\e[C_{i,n}(r)C_{i+m,n}(s)]\Big)+o(1)\\
&=& F(r)+\int_0^r \beta^{-1/2} (y) g^T(y)H^{-1}(y)\int_y^r \beta^{-1/2} (x) g(x)F(dx)F(dy)\\
&&+\int_0^r \beta^{-1/2} (y) g^T(y)H^{-1}(y)\int_y^s \beta^{-1/2} (x) g(x)F(dx)F(dy)\\
&&+ \int_0^r\int_0^s \beta^{-1/2} (y_1) g^T(y_1)H^{-1}(y_1)H(y_1\vee y_2)H^{-1}(y_2) \\
&&\hspace{10pt}\times\beta^{-1/2} (y_2)g(y_2)F(dy_2)F(dy_1)
+o(1) \: ,\\
&=& F(r)+o(1) \: ,
\end{eqnarray*}
which implies the assertion of the theorem. \hfill $\Box$

\section{The martingale transform of the process $\{ \Lambda_n(t) \}_{t \in [0,1]}$}
\def\theequation{4.\arabic{equation}}
\setcounter{equation}{0}

As pointed out in Section 2, the process $\{\sqrt{n}(\hat S_t (1 / \beta)-S_t (1 / \beta)) \}_{t \in [0,1]}$ (or its asymptotically equivalent
counterpart $\{ A^0_n (t) \}_{t \in [0,1]}$) depends on the unknown function $\beta$ (more precisely on the (unknown) functions
$\sigma^2(\cdot)$ and $m_4(\cdot)$). Similarly, the operator $T$ defined by (\ref{3.3}) is not completely known and has to be estimated from
the data. In this section we propose an empirical process, where the unknown quantities have been replaced by estimates and study the
application of an empirical version of the martingale transform. For this purpose we first have to specify the estimate in the process $\{
\Lambda_n(t)\}_{t \in [0,1]}$ defined in (\ref{proc}). We consider the Nadaraya-Watson weights
\begin{equation}
\label{3.11}w_{ij}=\frac{K\lt(\frac{t_{j,n}-t_{i,n}}{h}\rt)}{\sum_{l=1}^nK\lt(\frac{t_{l,n}-t_{i,n}}{h}\rt)} \: , \quad \quad i,j=1,\dots,n,
\end{equation}
at the points $t_{i,n} \ (i = 1, \dots, n)$ where $K$ denotes a symmetric kernel function and $h$ defines a bandwidth converging to 0 with
increasing sample size. The estimate of the function $\beta (\cdot)$ is now defined by
\begin{eqnarray} \nonumber
\hat{\beta}(t_{i,n})&=&\sum_{j=1}^nw_{ij}(Y_{j,n}-\hat{m}_h(t_{j,n}))^4 \\
&& + (4 \delta_r - 1) \sum^{n-r-1}_{j=1} w_{ij} (Y_{j,n} - \hat m_h
(t_{j,n}))^2 (Y_{j+r+1,n}- \hat m_h (t_{j+r+1,n}))^2, \label{3.13}
\end{eqnarray}
where $\hat m_h (t_{i,n}) = \sum^n_{j=1} w_{ij} Y_{j,n}$ denotes the Nadaraya-Watson estimate at the point $t_{i,n} \ (i=1,\dots,n)$.
Throughout this paper we assume that
\begin{description}
\item [(H)] The bandwidth $h$ satisfies $h=h_n=O\big(n^{-\frac{1}{2\gamma+1}}\big)$, where $\gamma > \frac {1}{2}$ denotes the Lipschitz
constant defined in Section 2.
\item[(K)] The kernel $K$ is symmetric, nonnegative, supported on the interval $[-1,1]$ and
satisfies $K(u)\leq 1$ for all $u\in[-1,1]$ and $K(u)\geq \kappa $ for all $|u|\leq 1/2$, where $\kappa>0.$
\end{description}

It will be proved in the appendix that under these additional assumptions
\begin{equation}
\sup_{t \in [0,1]} \Big| \label{3.17a}
 \frac{1}{\sqrt{n}}\sum_{i=1}^n1_{\{t_{i,n}\leq t\}}Z_{i,n}\{\hat{\beta}(t_{i,n})-\beta(t_{i,n})\}
 \Big| = o_p (1) \: ,
\end{equation}
and similar arguments as given in Dette and Hetzler (2008) show that
\be \label{proc1}
\Lambda_n(t)=\sqrt{n}\left ( \hat S_t (1/ \hat \beta) - S_t
(1/\beta)\right ) = A^1_n (t) + o_p \left(1\right)
\ee
uniformly with respect to $t \in [0,1]$. In this representation, $\{ A^1_n (t)\}_{t \in
[0,1]}$ denotes the process obtained from $\{ A^0_n (t) \}_{t \in [0,1]}$ by replacing $\beta (t)$ by its estimate $\hat \beta (t)$ defined in
(\ref{3.13}) and the vector $B_t(1/\beta)$ and the matrix $A$ by their estimates $\hat{B}_t(1/\hat{\beta})$ and $\hat A$ defined in (\ref{2.4}) and (\ref{2.6}), that is
\begin{eqnarray}
A_n^1(t)&=&C_n^1(t)-D_n^1(t),\label{3.13a}
\end{eqnarray}
where
\begin{eqnarray}
C_n^1(t) &=&\frac{\sqrt{n}}{n-r}\sum_{i=r+1}^n1_{\{t_{i,n}\leq t\}}Z_{i,n}\hat{\beta}^{-1/2}(t_{i,n}) \label{2.18a} \: , \\
\nonumber D_n^1(t) &=& \hat{B}_t^T (1/\hat \beta) \hat{A}^{-1}\frac{\sqrt{n}}{n-r}\sum_{i=r+1}^{n}Z_{i,n}g(t_{i,n}) \: .
\end{eqnarray}
Similarly, we replace the operator $T$ by its empirical version defined by
\be \label{3.14}
(T_n\eta)(t)=\eta(t)-\int_0^t \hat \beta^{-1/2} (y)
g^T(y)H_n^{-1}(y)\int_y^1 \hat \beta^{-1/2} (z) g(z)\eta(dz)F_n(dy),
\ee where the matrix $H_n(x)$ is given by \be \label{3.15}
H_n(x)=\int_x^1 \hat \beta^{-1} (u) g(u)g^T(u)F_n(du)=\frac{1}{n}\sum_{i=1}^n1_{\{t_{i,n}\geq x\}} \hat \beta^{-1} (t_{i,n}) g(t_{i,n})g^T(t_{i,n}) \ee
 and $F_n (t) = \frac {1}{n} \sum^n_{i=1} 1_{\{ t_{i,n} \leq t \}}$ denotes the empirical distribution function of the
design points.

Note that the matrix $H(x)$ used in the transformation (\ref{3.3}) is singular at the point $x=1$, and as a consequence, the matrices
$H_n^{-1}(x)$ are unbounded on the whole interval $[0,1]$. To circumvent this difficulty, we restrict the process $T_nA_n^1$ to the interval
$[0,t_0]$ with a fixed $0<t_0<1.$ This approach was also suggested by Khmaladze (1993) and Stute, Thies and Zhu (1998) among others.

The following results show that the asymptotic properties of the processes $\{T A^0_n (t)\}_{t \in [0,t_0]}$ and \linebreak $\{T_n
A^1_n (t)\}_{t \in [0,t_0]}$ coincide, and as a consequence we obtain weak convergence of the martingale transform of the process defined on the left hand side of
(\ref{proc1}).

%Note that the matrix $H(x)$ used in the transformation (\ref{3.3}) is singular at the point $x=1$, and as a consequence, the matrices
%$H_n^{-1}(x)$ are unbounded on the whole interval $[0,1]$. To circumvent this difficulty, we restrict the process $T_nA_n^1$ to the interval
%$[0,t_0]$ with a fixed $0<t_0<1.$ This approach was also suggested by Khmaladze (1989) and Stute, Thies and Zhu (1998) among others.

\bigskip

{\bf Theorem 4.1.} {\it If the assumptions stated at the beginning of Section 2 and the assumptions (H) and (K) are satisfied, then for any $0
< t_0 < 1$ the process $\{ T_n A^1_n (t)\}_{t \in [0,t_0]}$ converges weakly on $D[0, t_0]$ to a Brownian motion in time $F$, that is}
\[\{T_nA_n^1(t)\}_{t \in [0,t_0]}\Dkonv \{W\circ F(t)\}_{t \in [0,t_0]}\]

\bigskip

{\bf Corollary 4.2.} {\it If the assumptions of Theorem 4.1 are satisfied, then for any $0 < t_0 < 1$ the process $\{ T_n \Lambda_n(t) \}_{t
\in [0,t_0]}$ converges weakly on $D [0, t_0]$ to a Brownian motion in time $F$, that is}
\[\{ T_n \Lambda_n(t) \}_{t \in [0,t_0]} \Dkonv \{ W \circ F (t) \}_{t \in [0,t_0]}.\]
\bigskip

{\bf Proof of Theorem 4.1.} Obviously the assertion follows from the statement \be \label{3.16}
\sup_{t\in[0,t_0]}|TA_n^0(t)-T_nA_n^1(t)|=o_p(1) \: . \ee In order to prove the estimate (\ref{3.16}) we note that (using the notation $D^1_n =
A^1_n - C^1_n$)
\begin{eqnarray*} (T_nC_n^1-T_nA_n^1) (t)&=&
 D_n^1 (t) -\int_0^t \hat \beta^{-1/2}(y) g^T(y)H_n^{-1}(y)\int_y^1 \hat \beta^{-1/2}(z) g(z)D_n^1(dz)F_n(dy)\\
&=& D_n^1 (t) -\int_0^t \hat \beta^{-1/2}(y) g^T(y)F_n(dy)\hat{A}^{-1}\nonumber\big(\frac{\sqrt{n}}{n-r}\sum_{i=r+1}^nZ_{i,n}g(t_{i,n})\big)\nonumber\\
&=& D_n^1 (t) -\hat{B}_t (1/\hat \beta) \hat{A}^{-1}\frac{\sqrt{n}}{n-r}\sum_{i=r+1}^nZ_{i,n}g(t_{i,n})\nonumber =0.
\nonumber\label{tn}
\end{eqnarray*}
Consequently (observing the corresponding result for $T A^0_n - T C^0_n$ in (\ref{z1})),
the assertion follows if the statement \be \label{3.17} \sup_{t\in[0,t_0]}|TC_n^0(t)-T_nC_n^1(t)|=o_p(1) \ee can be proved,
where
\begin{eqnarray*}
TC_n^0(t)&=&C_n^0(t)-\int_0^t \beta^{-1/2}(y) g^T(y)H^{-1}(y)\int_y^1 \beta^{-1/2}(z) g(z)C_n^0(dz)F(dy)
= C_n^0(t)-B_n^0(t), \\
T_nC_n^1(t)&=&C_n^1(t)-\int_0^t \hat \beta^{-1/2}(y) g^T(y)H_n^{-1}(y)\int_y^1 \hat \beta^{-1/2}(z) g(z)C_n^1(dz)F_n(dy)=
 C_n^1(t)-B_n^1(t),
 \end{eqnarray*}
$C^0_n$ and $C^1_n$ are defined in (\ref{2.17}) and (\ref{2.18a}), respectively, and the equalities define the processes $B^0_n$ and $B^1_n$ in an obvious manner. It follows by a Taylor expansion, by the estimate (\ref{3.17a}) and the
estimate
\begin{equation} \label{mack}
 \sup_{t\in[0,t_0]}|\hat{m}_h(t)-m(t)|=O_p\lt(n^{-\frac{\gamma}{2\gamma+1}}\sqrt{\log n}\rt)
 \end{equation}
[see Mack and Silverman (1982)] that
\be \label{3.18} \sup_{t\in[0,t_0]}|C_n^1(t)-C_n^0(t)|=o_p(1). \ee
We now consider the remaining
difference
\begin{eqnarray} \nonumber
{B_n^1(t)-B_n^0(t)}
&=& \int_0^t\hat{\beta}^{-1/2}(y)g^T(y)H_n^{-1}(y)\hat{U}_n(y)F_n(dy)-\int_0^t\beta^{-1/2}(y)g^T(y)H^{-1}(y)U_n(y)F(dy)\nonumber\\
&=&\int_0^t\big(\hat{\beta}^{-1/2}(y)-\beta^{-1/2}(y)\big)g^T(y)(H_n^{-1}(y)-H^{-1}(y))(\hat{U}_n(y)-U_n(y))F_n(dy)\nonumber\\
&&+ \int_0^t\beta^{-1/2}(y)g^T(y)(H_n^{-1}(y)-H^{-1}(y))(\hat{U}_n(y)-U_n(y))F_n(dy)\nonumber\\
&&+ \int_0^t\big(\hat{\beta}^{-1/2}(y)-\beta^{-1/2}(y)\big)g^T(y)H^{-1}(y)(\hat{U}_n(y)-U_n(y))F_n(dy)\nonumber\\
&&+ \int_0^t\beta^{-1/2}(y)g^T(y)H^{-1}(y)(\hat{U}_n(y)-U_n(y))F_n(dy)\nonumber\\
&&+\int_0^t\big(\hat{\beta}^{-1/2}(y)-\beta^{-1/2}(y)\big)g^T(y)(H_n^{-1}(y)-H^{-1}(y))U_n(y)F_n(dy)\nonumber\\
&&+\int_0^t\beta^{-1/2}(y)g^T(y)(H_n^{-1}(y)-H^{-1}(y))U_n(y)F_n(dy)\nonumber\\
&&+\int_0^t\big(\hat{\beta}^{-1/2}(y)-\beta^{-1/2}(y)\big)g^T(y)H^{-1}(y)U_n(y)F_n(dy)\nonumber\\
&&+ \int_0^t\beta^{-1/2}(y)g^T(y)H^{-1}(y)U_n(y)F_n(dy)-\int_0^t\beta^{-1/2}(y)g^T(y)H^{-1}(y)U_n(y)F(dy)\nonumber\\
&=& T_{n,1}(t)+\ldots+T_{n,7}(t) \label{3.19}\\
&&+ \int_0^t\beta^{-1/2}(y)g^T(y)H^{-1}(y)U_n(y)F_n(dy)-\int_0^t\beta^{-1/2}(y)g^T(y)H^{-1}(y)U_n(y)F(dy),\nonumber
\end{eqnarray}
where the last equality defines the terms $T_{n,1}(t), \dots, T_{n,7}(t)$ in an obvious manner and we have used the notation
$$
\hat{U}_n(y)= \int_y^1 \hat \beta^{- 1/2}(z) g(z)C_n^1(dz) \: , \qquad
U_n(y)= \int_y^1 \beta^{-1/2}(z) g(z)C_n^0(dz) \: .
$$
The nine terms in this expression are estimated separately. Because the proceeding for the first seven terms is similar, we exemplarily illustrate the arguments for $T_{n,6}(t)$. This term is bounded by
\[\sup_{y\in[0,t_0]}\|H_n^{-1}(y)-H^{-1}(y)\|T_{n1},\]
where
\begin{eqnarray*} T_{n1} &:=& \int_0^t | \beta^{-1/2} (y) | \: \| g^T(y)\| \| U_n(y)\| F_n(dy) \\
& \leq & T_{n11} := \frac{1}{n}\sum_{i=1}^n \Biggl[ | \beta^{-1/2} (t_{i,n}) | \: \big\| g^T(t_{i,n})\big\| \big\|
\frac{\sqrt{n}}{n-r}\sum_{j=r+1}^n1_{\{t_{j,n}\geq t_{i,n}\}}g(t_{j,n})Z_{j,n} \beta^{-1}(t_{j,n})\big\| \Biggr] ,
 \end{eqnarray*}
$\parallel\cdot\parallel$ denotes the euclidean norm on $\mathbb{R}^d$ and its induced matrix norm on $\mathbb{R}^{d \times d}$
simultaneously, and we have used the definition of $U_n (y)$. A straightforward application of the Cauchy-Schwarz inequality shows
\begin{eqnarray*}
\e T_{n11} &\leq & \frac{1}{n}\sum_{i=1}^n | \beta^{-1/2} (t_{i,n}) | \: \big\| g^T(t_{i,n})\big\|\bigg(\e\|
\frac{\sqrt{n}}{n-r}\sum_{j=r+1}^n1_{\{t_{j,n}\geq t_{i,n}\}}g(t_{j,n})Z_{j,n} \beta^{-1}(t_{j,n})\|^2\bigg)^{1/2} = O(1)
\end{eqnarray*}
uniformly with respect to $t \in [0, t_0]$, which implies $T_{n1} = O_p (1)$ uniformly on the interval $[0, t_0]$. Using similar arguments as in the proof of the estimate (\ref{3.17a}) it can be shown that
\[\max_{i=1,\ldots,n}|\hat{\beta}(t_{i,n})-\beta(t_{i,n})|=o_p(1).\]
By a Taylor expansion and the assumption (\ref{1.2}) on the design it follows that \be\label{3.20} \sup_{y\in[0,t_0]}\|H_n^{-1}(y)-H^{-1}(y)\| = o (1), \ee and we obtain that
$T_{n,6}(t)$ is of order $o_p (1)$ uniformly with respect to $t\in[0,t_0]$. Using similar arguments it follows that $T_{n,1}(t),
\ldots,T_{n,5}(t)$ and $T_{n,7}(t)$ are also of order $o_p (1)$ uniformly in $t\in[0,t_0]$. For the difference of the last two terms on the
right-hand side of (\ref{3.19}) we show the estimate
\be\label{3.21}
\sup_{t\in[0,t_0]}\bigg|\int_0^t \beta^{-1/2} (y)
g^T(y)H^{-1}(y)U_n(y)\bigg( F_n(dy)-F(dy) \bigg) \bigg|=o_p(1)
\ee
 using Lemma 6.6.4 in Koul (2002). Note that for the application of this result one has to show the tightness of the process $\{ U_n (x)
\}_{x \in [0, t_0]}$. For this purpose we consider the components of $U_n$ separately, that is
\[U_n^{(p)}(x)=\frac{\sqrt{n}}{n-r}\sum_{i=r+1}^n1_{\{t_{i,n}\geq x\}}\sigma_p^2(t_{i,n})\beta^{-1}(t_{i,n})Z_{i,n}
,~\quad p=1,\ldots,d,\] and introduce the notation
\[\nu_p(x)=1_{\{y_1\leq x< y_2\}}\sigma_p^2(x)\beta^{-1}(x).\]
Now a similar calculation as in Dette and Hetzler (2008) yields
\begin{eqnarray*}\e[(U_n^{(p)}(y_2)-U_n^{(p)}(y_1))^4]
&=&\frac{n^2}{\left(n-r\right)^4}\e\Big[\Big(\sum_{i=r+1}^n\nu_p(t_{i,n})Z_{i,n}\Big)^4\Big] \leq C(y_2-y_1)^2
\end{eqnarray*}

for some constant $C>0$ and $0\leq y_1\leq y_2 \leq t_0$. This implies tightness of each component $U_n^{(p)}$ [see Billingsley (1999)] and as
a consequence tightness of the process $U_n$ [see Billingsley (1979)]. \hfill $\Box$

\bigskip

{\bf Remark 4.3.} Theorem 4.1 and Corollary 4.2 remain correct if the Nadaraya-Watson weights in the estimate $\hat \beta$ defined in (\ref{3.13}) are replaced by local
linear weights. This follows by a careful inspection of the proof of the estimate (\ref{3.17a}) in the appendix. In practical applications the
use of local linear weights is strictly recommended because of the better performance of the local linear estimate at the boundary of the
design space.

\bigskip

\section{Finite sample properties}
\def\theequation{5.\arabic{equation}}
\setcounter{equation}{0}

In this section we investigate the finite sample properties of the new test by means of a simulation study. We have generated data according to
the model
\be \label{sim0}
Y_{i,n} = 1 + t_{i,n} + \sigma (t_{i,n}) \eps_{i,n},\quad i=1,\dots,n,
\ee
where $t_{i,n} = i/(n+1),~~i=1, \dots,n$, and simulated the power of the test for the hypothesis
\be \label{sim1}
H_0:\ \sigma^2 (t) = 1 + \theta t^2 \
\ee
and the variance functions
\begin{eqnarray}
 \sigma^2(t)&=& 0.5 + 3 t^2 + 2.5c\sin(2\pi t),\label{sim2}\\
 \sigma^2(t)&=& 0.5 + 3 t^2 + 2ce^{2t},\label{sim3}\\
 \sigma^2(t)&=& 0.5 + 3 t^2 + 4c\sqrt{t}.\label{sim4}
\end{eqnarray}
Note that the choice $c=0$ in (\ref{sim2}) - (\ref{sim4}) corresponds to the null hypothesis of a quadratic variance function. The errors
$\eps_{i,n}$ are standard normal distributed and we use a difference sequence of order $r=1$ for the calculation of the pseudo residuals
$R_{i,n}$, which determines the weights as $d_0=-d_1=1/\sqrt{2}$ and yields $\beta(x)=m_4(x)\sigma^4(x).$ In order to apply the test we have
to calculate the transformation
\[T_n\Lambda_n(t)=T_n\big(\sqrt{n}\big(\hat{S}_t(1/\hat{\beta})-S_t(1/\beta)\big)\big)\]
for the process $\Lambda_n(t)$ given in (\ref{proc1}). Under the null hypothesis (\ref{sim1}) we have $S_t(1/\beta)=0$ for
all $t\in[0,1],$ and the process $\Lambda_n(t)$ can be written as
\[\Lambda_n(t)=\sqrt{n}\,\hat{S}_t(1/\hat{\beta})=\hat C_n(t)-\hat D_n(t) \]
with
\begin{eqnarray*}
\hat C_n(t)&=&\frac{\sqrt{n}}{n-r}\sum_{i=r+1}^n1_{\{t_{i,n}\leq t\}}\hat{\beta}^{-1/2}(t_{i,n})R_{i,n}^2,\\
\hat D_n(t)&=&\hat B_t^T(1/\hat{\beta})\hat A^{-1}\frac{\sqrt{n}}{n-r}\sum_{i=r+1}^{n}R_{i,n}^2g(t_{i,n}).
\end{eqnarray*}
By a similar argument as given in the proof of Theorem 4.1 it can be shown that $T_n\hat D_n(t)=0$ for all $t\in[0,1]$, and as a consequence
it is sufficient to calculate the transformation $T_n\hat C_n$. We use the Cram\'er-von-Mises statistic
$G_n=\int_0^1(T_n\Lambda_n)^2(t)dF_n(t),$ and from Corollary 4.2 and the Continuous Mapping Theorem it follows that \be \label{sim5}
G_n=\int_0^1(T_n\Lambda_n)^2(t)\,dF_n(t)\Dkonv \int_0^1W^2(F(t))\,dF(t)=\int_0^1W^2(t)\,dt, \ee where $W$ denotes a standard Brownian motion.
If $w_{\alpha}$ denotes the $1-\alpha$ quantile of the distribution of the random variable $\int^1_0 W^2(t)dt$, then the test, which rejects
the null hypothesis (\ref{sim1}) if
\begin{equation}\label{test}
G_n\geq w_{\alpha}
\end{equation}
 has asymptotically level $\alpha$ and is consistent against local alternatives converging to the null hypothesis at a rate $n^{-1/2}$.
As an estimator of the function $\beta(x)=m_4(x)\sigma^4(x)$ we use the estimator (\ref{3.13}), where $\hat m_h(\cdot)$ is the local
linear estimator of the regression function. The bandwidth for the calculation of the local linear estimate was determined
by least squares cross validation. If $h_{CV}$ is the bandwidth obtained by this procedure, the bandwidth in the
estimator (\ref{3.13}) was chosen as $h_{CV}/2$.

1000 simulation runs were performed in each scenario to calculate the rejection probabilities, which are shown in Table 5.1. For the sake of
comparison, the table also contains the corresponding rejection probabilities of the bootstrap test proposed by Dette and Hetzler (2008), which
are displayed in brackets. If the null hypothesis is satisfied $(c=0)$, we observe a rather precise approximation of the nominal level in all cases, even for the sample size $n=50$. Under
the alternatives the behaviour of the two tests is different. In model (\ref{sim2}) the bootstrap test yields a substantially larger power than
the test based on the martingale transformation, in particular for the sample size $n=50$ or $n=100$. In model (\ref{sim4}) the situation is similar for the sample size $n=50$, but the differences
between the rejection probabilities of the two steps are smaller. Moreover, for the sample size $n=200$, the test based on the martingale transform
shows a better performance. In model (\ref{sim3}) the bootstrap test is more powerful for the sample size $n=50$, while for the 
sample sizes $n=100$ and $n=200$ the test based on the martingale transformation always yields a larger power than the bootstrap test.

{\small
\begin{center}
\begin{tabular}{|c|c|r|r|r|r|r|r|r|r|r|}
\hline
&&\multicolumn{3}{|c|}{$n=50$} &\multicolumn{3}{|c|}{$n = 100$} &\multicolumn{3}{|c|}{$n = 200$}\\
\cline{3-11}
&c&.025 &.05 &.10&.025 &.05 &.10 &.025 &.05 &.10\\
\hline
&0& .027& .047& .101& .020& .042& .090& .022& .034& .078\\
&&(.030)&(.053)&(.105)&(.018)&(.044)&(.099)&(.024)&(.047)&(.096)\\
(\ref{sim2})&0.5 & .138& .230& .378& .313& .471& .631& .732 & .844 & .920\\
& & (.476)& (.549)&(.623)& (.702)& (.764)& (.820)& (.885) & (.923) & (.961)\\
&1& .258& .380& .559& .608& .743& .879& .958& .987 & .999\\
& & (.844)& (.874)& (.907)& (.976)& (.982)& (.990)& (.999) & (1.00) & (1.00)\\
\hline\hline
&0& .026& .045& .090& .026& .049& .086& .023& .047& .090\\
&&(.028)&(.049)&(.093)&(.034)&(.059)&(.106)&(.023)&(.046)&(.086)\\
(\ref{sim3})&0.5 & .097& .175& .292& .246& .359& .509& .520 & .646 & .785\\
& & (.193)& (.264)&(.372)& (.266)& (.346)& (.457)& (.398) & (.498) & (.619)\\
&1& .147& .235& .362& .339& .447& .610& .745& .847 & .936\\
& & (.243)& (.330)& (.426)& (.331)& (.418)& (.524)& (.561) & (.654) & (.747)\\
\hline \hline
&0& .029& .056& .107& .027& .047& .087& .027& .048& .094\\
&&(.019)&(.041)&(.097)&(.025)&(.044)&(.098)&(.023)&(.047)&(.102)\\
 (\ref{sim4})&0.5 & .088& .159& .260& .239& .355& .473& .517 & .646 & .765\\
& & (.231)& (.303)&(.400)& (.327)& (.429)& (.556)& (.537) & (.637) & (.734)\\
&1& .162& .273& .414& .386& .532& .685& .825& .910 & .957\\
& & (.381)& (.484)& (.586)& (.565)& (.667)& (.756)& (.775) & (.847) & (.911)\\
\hline
\end{tabular}
\end{center}}

{\bf Table 5.1.} {\it Rejection probabilities of the Cram\'{e}r-von-Mises test (\ref{test}) for the hypothesis (\ref{sim1}) in the regression
model (\ref{sim0}). The corresponding rejection probabilities of the bootstrap test proposed by Dette and Hetzler (2008) are displayed in
brackets.}

\vskip .75cm

{\bf Acknowledgements.} The work of the authors
 was supported by the Deutsche Forschungsgmeinschaft (SFB 475,
Komplexit\"atsreduktion in multivariaten
 Datenstrukturen, Teilprojekt B1)
 and in part by a NIH grant award
IR01GM072876:01A1. The authors are also grateful
%to a referee and the associate editor for their constructive
% comments on an earlier version of this paper and
to Martina Stein, who typed parts of this paper with considerable expertise.

\section{Appendix: Proof of (\ref{3.17a})}
\def\theequation{6.\arabic{equation}}
\setcounter{equation}{0}
Throughout this section
 we omit the index $n$; in particular we write $t_j$ and $Z_j$ instead of $t_{j,n}$ and $Z_{j,n},$ respectively.
For the sake of brevity we only indicate the main steps of the proof, details can be found in Hetzler (2008).
Furthermore we restrict ourselves
to the case $\sigma\equiv 1$ and $r=1$ and note
 that the general case is proved exactly in the same way with
 some additional notation.
 This simplification yields for the random variables $Z_i$
\[Z_i=d_0\sigma(t_i)\eps_i+d_1\sigma(t_{i-1})\eps_{i-1}=\frac{\eps_i-\eps_{i-1}}{\sqrt{2}}.\]
A straightforward calculation gives
\be
\label{a1}
A(t) = \frac{1}{\sqrt{n}}\sum_{i=1}^n1_{\{t_i\leq t\}}Z_i\{\hat{\beta}(t_i)-\beta(t_i)\}
=\sum_{j=1}^5A_j(t),
\ee
where
\begin{eqnarray*}
A_1(t)&=& \frac{1}{\sqrt{n}}\sum_{i=1}^n1_{\{t_i\leq t\}}Z_i\sum_{j=1}^nw_{ij}(\eps_j^4-m_4(t_i)),\\
A_2(t)&=& \frac{4}{\sqrt{n}}\sum_{i=1}^n1_{\{t_i\leq t\}}Z_i\sum_{j=1}^nw_{ij}\eps_j^3(m(t_j)-\hat{m}_h(t_j)),\\
A_3(t)&=& \frac{6}{\sqrt{n}}\sum_{i=1}^n1_{\{t_i\leq t\}}Z_i\sum_{j=1}^nw_{ij}\eps_j^2(m(t_j)-\hat{m}_h(t_j))^2,\\
A_4(t)&=& \frac{4}{\sqrt{n}}\sum_{i=1}^n1_{\{t_i\leq t\}}Z_i\sum_{j=1}^nw_{ij}\eps_j(m(t_j)-\hat{m}_h(t_j))^3,\\
A_5(t) &=& \frac{1}{\sqrt{n}}\sum_{i=1}^n1_{\{t_i\leq t\}}Z_i\sum_{j=1}^nw_{ij}(m(t_j)-\hat{m}_h(t_j))^4.
\end{eqnarray*}
We rewrite
$ m(t_j)-\hat{m}_h(t_j) = \rho_j-\sum_{k=1}^nw_{jk}\eps_k $
with \be\label{rho} \rho_j:=m(t_j)-\sum_{k=1}^n w_{jk}m(t_k)=\sum_{k=1}^n w_{jk}\big(m(t_j)-m(t_k)\big) \ee and
 first consider the term $A_1$. For its expectation we have
$$
\e A_1(t)
= \frac{1}{\sqrt{n}}\sum_{i=1}^n1_{\{t_i\leq t\}}\e[Z_i\sum_{j=1}^nw_{ij}h_{ij}],
$$
where we used the notation $h_{ij}:=\eps_j^4-m_4(t_i).$ Note that
$ |\e h_{ij}|=|m_4(t_j)-m_4(t_i)|\leq Lh^{\gamma}$
whenever $|t_j-t_i|\leq h$ (recall the H\"older continuity for the function $m_4$) and that
 it follows from the assumption on the design and the kernel
\begin{equation}
\label{a5}\frac{K_h\lt(t_j-t_i\rt)}{C_2n}\leq w_{ij}\leq \frac{K_h\lt(t_j-t_i\rt)}{\kappa C_1n/2}
\end{equation}
where $K_h(x) = K(x/h)/h$ and $C_1$ and $C_2$ denote positive constants. This yields
\[\e[Z_i\sum_{j=1}^nw_{ij}h_{ij}]=\e[Z_i(w_{ii}h_{ii}+w_{i,i-1}h_{i,i-1})]=O\lt(\frac{1}{nh}\rt)\]
uniformly with respect to $i=r+1,\ldots,n$, and it follows that
$\e A_1(t)=O(\frac{1}{h\sqrt{n}})=o(1).$
For the estimation of the second moment we decompose $A_1^2(t)$ as follows
\begin{eqnarray*}
A_1^2(t) &=& D_1(t)+D_2(t)+D_3(t),
\end{eqnarray*}
with
\begin{eqnarray*}
D_1(t) &=& \frac{1}{n}\sum_{i=1}^n1_{\{t_i\leq t\}}Z_i^2(\sum_{j=1}^nw_{ij}h_{ij})^2,\\
D_2(t) &=& \frac{2}{n}\sum_{i=1}^{n-1}1_{\{t_i\leq t\}}1_{\{t_{i+1}\leq t\}}Z_iZ_{i+1}\sum_{j=1}^nw_{ij}h_{ij}\sum_{k=1}^nw_{i+1,k}h_{i+1,k},\\
D_3(t) &=& \frac{1}{n}\sum_{|i-l|\geq 2}1_{\{t_i\leq t\}}1_{\{t_l\leq t\}}Z_iZ_l\sum_{j=1}^nw_{ij}h_{ij}\sum_{k=1}^nw_{lk}h_{lk}.\\
\end{eqnarray*}
Observing (\ref{a5}) it follows for the set $A_i:=\{i-1,i\}$
\begin{eqnarray*}
\e[Z_i^2\sum_{j,k=1}^nw_{ij}w_{ik}h_{ij}h_{ik}]
&=&\e[Z_i^2\sum_{j,k\n A_i}w_{ij}w_{ik}h_{ij}h_{ik}]+ 2\e[Z_i^2(w_{ii}h_{ii}+w_{i,i-1}h_{i,i-1})\sum_{k\n A_i}w_{ik}h_{ik}]\\
&&+ \e[Z_i^2(w_{ii}h_{ii}+w_{i,i-1}h_{i,i-1})^2] \\
&=& \e[Z_i^2\sum_{j,k\n A_i}w_{ij}w_{ik}h_{ij}h_{ik}] + O\lt(n^{-1}h^{-1}\rt) + O\lt(n^{-2}h^{-2}\rt) \\
%&&\\
&=& O(h^{2\gamma}) + O\lt(n^{-1}h^{-1}\rt) + O\lt(n^{-2}h^{-2}\rt) ~\\
&=& O(h^{2\gamma}).
\end{eqnarray*}
A similar calculation shows
$ ED_2(t)=O(h^{2\gamma}). $
For the remaining estimate for the term $D_3(t)$ we consider the set $A_{i,l}=\{i-1,i,l-1,l\}$ and obtain
\begin{eqnarray*}
{ \e [Z_iZ_l\sum_{j=1}^nw_{ij}h_{ij}\sum_{k=1}^nw_{lk}h_{lk}]}
&=&\e[Z_iZ_l\sum_{j,k\n A_{i,l}}w_{ij}w_{lk}h_{ij}h_{lk}]\\
&&+ \e[Z_iZ_l(w_{ii}h_{ii}+w_{i,i-1}h_{i,i-1}+w_{il}h_{il}+w_{i,l-1}h_{i,l-1})\sum_{k\n A_{i,l}}w_{lk}h_{lk}]\\
&&+ \e[Z_iZ_l(w_{li}h_{li}+w_{l,i-1}h_{l,i-1}+w_{ll}h_{ll}+w_{l,l-1}h_{l,l-1})\sum_{j\n A_{i,l}}w_{ij}h_{ij}]\\
&&+\e[Z_iZ_l\sum_{j,k\in A_{i,l}}w_{ij}w_{lk}h_{ij}h_{lk}].
\end{eqnarray*}
Note that the random variables
 $Z_i$ and $Z_l$ are independent whenever $|l-i|\geq 2$ and consequently the first three terms
 in the above expression vanish. The remaining fourth term can be decomposed in a sum of $16$, which are
 all of the form
\bea \e[Z_iZ_lw_{ii}w_{ll}h_{ii}h_{ll}]=O\lt(\frac{1}{n^2h^2}\rt). \eea This yields
\[\e D_3(t)=O\lt(\frac{1}{nh^2}\rt)\]
and as a consequence $\e A_1^2(t)=O\lt(\frac{1}{nh^2}\rt)=o(1)$. Thus we obtain
\begin{equation} \label{a01}
A_1(t)=o_p(1)
\end{equation}
uniformly with respect to $t\in[0,t_0].$ In order to derive a corresponding
estimate for the term $A_2$ we use the decomposition
\begin{eqnarray*}
A_2(t)&=& A_{21}(t)-A_{22}(t)
\end{eqnarray*}
with
\begin{eqnarray*}
 A_{21}(t) &=& \frac{4}{\sqrt{n}}\sum_{i=1}^n1_{\{t_i\leq t\}}Z_i\sum_{j=1}^nw_{ij}\eps_j^3\rho_j,\\
A_{22}(t) &=& \frac{4}{\sqrt{n}}\sum_{i=1}^n1_{\{t_i\leq t\}}Z_i\sum_{j=1}^nw_{ij}\eps_j^3\sum_{k=1}^nw_{jk} \eps_k.
\end{eqnarray*}
Now the H\"older continuity of the regression function implies
% \begin{equation} \label{a2}
$ |\rho_j|=|\sum_{k=1}^nw_{jk}(m(t_j)-m(t_k))|\leq Lh^\gamma$
%\end{equation}
for some positive constant $L$ and a straightforward calculation shows (note that the random variables $Z_i$ depend only on $\eps_i$ and $\eps_{i-1}$)
\[\e[Z_i\sum_{j=1}^nw_{ij}\eps_j^3\rho_j]=O\lt(\frac{h^{\gamma-1}}{n}\rt),\]
which implies
\begin{equation}
\label{a3} \e[A_{21}(t)]=O\lt(\frac{h^{\gamma-1}}{\sqrt{n}}\rt).
\end{equation}
By a similar calculation it follows that
% \begin{eqnarray*}
% \e[Z_i\sum_{j=1}^nw_{ij}\eps_j^3\sum_{k=1}^nw_{jk}\eps_k] &=&
% \e[Z_i\sum_{j=1}^nw_{ij}w_{jj}\eps_j^4]+\e[Z_i\sum_{k\not =j}w_{ij}w_{jk}\eps_j^3\eps_k] \\
% &=& O\lt(\frac{1}{n^2h^2}\rt) +\e[Z_i\sum_{k\not =j}w_{ij}w_{jk}\eps_j^3\eps_k] \\
% &=& O\lt(\frac{1}{n^2h^2}\rt) + O\lt(\frac{1}{nh}\rt) = O\lt(\frac{1}{nh}\rt),
% \end{eqnarray*}
% which implies
$\e[A_{22}(t)]=O\lt(\frac{1}{h\sqrt{n}}\rt) $ and a combination of this estimate with (\ref{a3}) gives
\begin{equation}
\label{a3a}\e[A_2(t)]=O\lt(\frac{1}{h\sqrt{n}}\rt).
\end{equation}
The estimation of the second moments of $A_{21}(t)$ and $A_{22}(t)$ is more complicated and we indicate the calculations for the term
$A_{21}(t)$, which can be decomposed as
\begin{eqnarray*}
A_{21}^2(t)&=& B_1(t)+B_2(t)+B_3(t),
\end{eqnarray*}
where
\begin{eqnarray*}
B_1(t) &=&
\frac{16}{n}\sum_{i=1}^n1_{\{t_i\leq t\}}Z_i^2\big(\sum_{j=1}^nw_{ij}\eps_j^3\rho_j\big)^2,\\
B_2(t) &=& \frac{32}{n}\sum_{i=1}^{n-1}1_{\{t_i\leq t\}}1_{\{t_{i+1}\leq t\}}Z_iZ_{i+1}\sum_{j=1}^nw_{ij}\eps_j^3\rho_j\sum_{l=1}^nw_{i+1,l}\eps_l^3\rho_l,\\
B_3(t) &=& \frac{16}{n}\sum_{|l-i|\geq 2}1_{\{t_i\leq t\}}1_{\{t_l\leq
t\}}Z_iZ_l\sum_{j=1}^nw_{ij}\eps_j^3\rho_j\sum_{r=1}^nw_{lr}\eps_r^3\rho_r .
\end{eqnarray*}
 Using the estimates (\ref{a5}) we obtain
\begin{eqnarray*}
E[B_1(t) ] &=&\frac{16}{n} \sum_{i=1}^n E \Bigl[ 1_{\{t_i\leq t\}}Z_i^2\sum_{j=1}^nw_{ij}^2\eps_j^6\rho_j^2 \Bigr] + \frac{16}{n} \sum_{i=1}^n
E \Bigl[ 1_{\{t_i\leq t\}}Z_i^2\sum_{j\ne l}w_{ij}w_{il}\eps_j^3\eps_l^3\rho_j\rho_l
 \Bigr] \\
&=& O\lt(\frac{h^{2\gamma-1}}{n}\rt) + O(h^{2\gamma}) ~= ~ O(h^{2\gamma}).
\end{eqnarray*}
A similar calculation shows $\e B_2(t)=O(h^{2\gamma})$ and
\[\e[B_3(t)]=O\lt(\frac{h^{2\gamma-2}}{n}\rt),\]
which implies
\begin{equation}\label{a6}
\e[A_{21}^2(t)]=O(h^{2\gamma}).
\end{equation}
Similarly we obtain
\[\e A_{22}^2(t)=O\lt(\frac{1}{nh^2}\rt),\]
and a combination with (\ref{a6}) gives
\[\e A_2^2(t)=O\lt(\frac{1}{nh^2}\rt)=O\lt(n^{\frac{1-2\gamma}{2\gamma+1}}\rt)=o(1). \]
On the other hand we have from (\ref{a3a}) the estimate $E A_2(t)=O\lt(\frac{1}{h\sqrt{n}}\rt)=O\lt(n^{\frac{1-2\gamma}{4\gamma+2}}\rt)=o(1)$
and it follows that
\begin{equation} \label{a02}
A_2(t)=o_p(1)
\end{equation}
uniformly on the interval $t\in[0,t_0].$
 The term $A_3(t)$ can be treated by similar arguments, which are omitted for the sake of brevity [see
Hetzler (2008) for more details]. Tedious calculations yield
\begin{equation} \label{a03}
A_3(t)=o_p(1)
\end{equation}
uniformly with respect to $t\in[0,t_0].$ Finally we use the estimate (\ref{mack})
 and obtain the remaining terms in (\ref{a1})
\begin{eqnarray*}
|A_4(t)|
&\leq & \frac{4}{\sqrt{n}}\sum_{i=1}^n\big|Z_i\sum_{j=1}^nw_{ij}\eps_j (m(t_j)-\hat{m}_h(t_j))^3\big|\\
&\leq & \sup_{t\in[0,t_0]}|\hat{m}_h(t)-m(t)|^3\cdot\frac{4}{\sqrt{n}}\sum_{i=1}^n|Z_i|\sum_{j=1}^nw_{ij}|\eps_j|\\
&=& O_p\lt(n^{-\frac{3\gamma}{2\gamma+1}}(\log{n})^{3/2}n^{\frac{2\gamma+2}{4\gamma+2}}\rt) ~=~
O_p\lt(n^{\frac{2-4\gamma}{4\gamma+2}}\log{n}\rt)~=~o_p(1)
\end{eqnarray*}
and
\begin{eqnarray*}
|A_5(t)|
&\leq & \sup_{t\in[0,t_0]}\lt|\hat{m}_h(t)-m(t)\rt|^4\frac{1}{\sqrt{n}}\sum_{i=1}^n |Z_i|\\
&=& O_p\lt(n^{-\frac{4\gamma}{2\gamma+1}}\sqrt{n}(\log{n})^2\rt)~ =~
 O_p\lt(n^{\frac{1-6\gamma}{4\gamma+2}}\log{n}\rt) ~=~ o_p(1)
\end{eqnarray*}
uniformly in $t\in [0,t_0]$. Combining these estimates with (\ref{a01}), (\ref{a02}) and (\ref{a03}) it follows that $A(t)=o_p(1)$ holds
uniformly with respect to $t\in [0,t_0]$, which proves (\ref{3.17a}).

\vskip 1cm

{\Large References}

\bigskip

N. I. Achieser (1956). Theory of Approximation. Frederik Ungar Publishing Co., N.Y.

\medskip

P. L. Bickel (1978). Using residuals robustly I: Tests for heteroscedasticity and nonlinearity. Annals of Statistics 6, 266-291.

\medskip

 {P. Billingsley (1979). Probability and measure. Wiley Series in Probability and Statistics, N.Y.}

\medskip

{P. Billingsley (1999). Convergence of probability measures, 2nd ed. Wiley Series in Probability and Statistics, N.Y.}

\medskip

T. S. Breusch and A. R. Pagan (1979). A simple test for heteroscedasticity and random coefficient variation. Econometrica 47, 1287-1294.

\medskip

R. D. Cook and S. Weisberg (1983). Diagnostics for heteroscedasticity in regression. Biometrika 70, 1-10.

\medskip

{H. Dette (2002). A consistent test for heteroscedasticity in nonparametric regression based on the kernel method. Jour.\ Statist.\ Plan.\ and
Infer.\ 103, 311-330.}

\medskip

Dette and Hetzler (2008). {A simple test for the parametric form of the variance function in nonparametric regression. To appear in: Annals of
the Institute of Statistical Mathematics. Online: http://www.springerlink.com/content/102845/?Content+Status=Accepted}

\medskip

{H. Dette and A. Munk (1998). Testing heteroscedasticity in nonparametric regression. J.\ R.\ Statist.\ Soc.\ B, 60, 693-708}

\medskip

H. Dette, I. van Keilegom and N. Neumeyer (2007). A new test for the parametric form of the variance function in nonparametric regression.
J.\ Roy.\ Statist.\ Soc., Ser.\ B 69, 903-917.

\medskip

{T. Gasser, L. Sroka and G. Jennen-Steinmetz (1986). Residual variance and residual pattern in nonlinear regression. Biometrika 73, 626-633.}

\medskip

{P. Hall, J.W. Kay and D.M. Titterington (1990). Asymptotically optimal difference-based estimation of variance in nonparametric regression.
Biometrika 77, 521-528.}

\medskip

B. Hetzler (2008). Tests auf parametrische Struktur der Varianzfunktion in der nichtparametri\-schen Regression. PhD Thesis, Fakult\"at f\"ur
Mathematik, Ruhr Universit\"at Bochum (in German).

%\medskip

%E.V. Khmaladze (1979). The use of $\omega^2$-tests for testing parametric hypothesis. Theory Probab.\ Appl.\ 24, 283-301.

\medskip

E.V. Khmaladze (1981). A martingale approach in the theory of goodness-of-fit tests. Theory Probab.\ Appl.\ 26, 240-257.

\medskip

E.V. Khmaladze (1993). Goodness of fit problem and scanning innovation martingales. Ann.\ Statist.\ 21, 798-829.

\medskip

E.V. Khmaladze and H.L. Koul (2004). Martingale transforms goodness-of-fit tests in regression models. Ann.\ Statist.\ 32, 995-1034.

\medskip

H.L. Koul (2002). Weighted empirical processes in dynamic nonlinear models. Lecture Notes in Statistics, 166. Springer, N.Y.
\medskip

H.L. Koul (2006). Model diagnostics via martingale transforms: a brief review. In: Frontiers in Statistics (ed. J. Fan, H.L. Koul), 183-206.
Hackensack, NJ. World Scientific.

\medskip

{H. Liero (2003). Testing homoscedasticity in nonparametric regression. J.\ Nonpar.\ Statist.\, 15, 31-51.}

\medskip

Y.P. Mack, B.W. Silverman (1982). Weak and strong uniform consistency of kernel regression estimates.
Z. Wahrsch. Verw. Gebiete 61, 405-415.

\medskip

{B. {\O}ksendal (2003). Stochastic Differential Equations: An Introduction with Applications. \linebreak Springer-Verlag.}

\medskip

A. R. Pagan and Y. Pak (1993). Testing for Heteroscedasticity in: G. S. Maddala, C. R. Rao and H. D. Vinod, eds., Handbook of Statistics, Vol.
11, Elsevier Science Publishers B.V.

\medskip

D. Pollard (1984). Convergence of Stochastic Processes. Springer, N.Y.

\medskip

{J. Sacks and D. Ylvisaker (1970). Design for regression problems with correlated errors III. Ann.\ Math.\ Statist.\ 41, 2057-2074.}

\medskip

W. Stute, S. Thies, L. Zhu (1998). Model checks for regression: an innovation process approach. Ann.\ Statist.\ 26, 1916-1934.

\medskip

{L. Zhu, Y. Fujikoshi and K. Naito (2001). Heteroscedasticity checks for regression models. Science in China (Series A), 44, 1236-1252}

\end{document}